

\baselineskip=14pt
\parskip=10pt

\font\eighttt=cmtt8
\magnification=\magstephalf

\def\1{{\overline{1}}}
\def\2{{\overline{2}}}
\parindent=0pt
\overfullrule=0in

\def\frac#1#2{{#1 \over #2}}


\centerline
{
\bf 
Zeroless Arithmetic: Representing Integers ONLY using ONE 
}
\bigskip
\centerline
{
\it Edinah K. GNANG and Doron ZEILBERGER
}

\quad \quad \quad \quad \quad
{\it
``The One counts Himself, and no-one else counts Him, and He is every number, He is root, and foundation and square and cube, and He is like the essence 
that carries all the cases, and every number is in His power, and He is in every number in deed, and He is present, and every number is present 
because of Him, and He is Ancient, and every [other] number is [re]newed, and He is the reason for every number, pair[even] and that is not pair, 
He is not a number, and will not multiply and will not divide.'' }
\medskip
\quad \quad \quad \quad \quad - Abraham Ibn Ezra (1089-1164), {\it Sefer HaEkhad} (``Book of One'')[I]

{\bf Abstract}: We use {\it recurrence equations} (alias difference equations) to
enumerate the number of formula-representations of positive integers using
only addition and multiplication, and using  addition, multiplication, and  exponentiation,
{\bf where all the inputs are ones}.
We also describe efficient algorithms for the random generation of such representations,
and use Dynamical Programming to find a shortest possible formula representing any given
positive integer.

{\bf Very Important}: This article is accompanied by the Maple package

{\eighttt http://www.math.rutgers.edu/\~{}zeilberg/tokhniot/ArithFormulas} \quad ,

and the  output files that are linked to from the webpage (``front'') of this article

{\tt http://www.math.rutgers.edu/\~{}zeilberg/mamarim/mamarimhtml/arif.html} \quad \quad .

{\bf Prologue}

According to conventional wisdom, the invention (``discovery'') of {\it zero} was one of the greatest moments
in the annals of mathematics. We respectfully disagree. The invention of zero was a great disaster, that lead
to the beginning of {\it nihilism}. Here we will show how it is possible to manage very well without $0$.

{\bf Introduction}

Mark Twain once wrote a letter to a friend that started with

{\it ``I didn't have time to write a short letter so I wrote a long one ...''}

We mathematicians (and computer scientists) deal with {\it numbers} rather than words, but even
the seemingly naive question of {\it representing} a positive integer as succinctly as possible
is far from trivial. 

This interesting question was addressed in [GD], where the systematic study of arithmetical formula-representation
was initiated, and two natural ways, called there ``the first canonical form'' (FCF), and the
``second canonical form'' (SCF) were introduced. The present article is a natural follow-up of [GD], but
in order to make it self-contained, we will review the basic notions.

We are all familiar with the ``caveman's representation'' of a positive integer by marking
lines (only using $1$'s), for example
$$
17=11111111111111111 \quad,
$$
also called the {\it unary} representation. More ``efficiently'' we have the familiar decimal,
`positional' systems that, alas, needs {\bf ten} symbols. 
The binary representation ``only'' uses two (one too many!) symbols, $0$ and $1$,
where,  for example, seventeen is written as $10001$, meaning
$$
17= 1 \cdot 2^4 +  0 \cdot 2^3 + 0 \cdot 2^2 + 0 \cdot 2^1 + 1 \cdot 2^0 \quad . 
$$
One can use the ``sparse notation'' by {\it only} keeping the 1's
$$
17=  2^4 + 1 \quad . 
$$
and doing the same for the exponents
$$
17=2^{2^2} + 1 \quad ,
$$
and finally replacing $2$ by $1+1$ getting an expression that {\bf only} uses $1$
$$
17=(1+1)^{(1+1)^{1+1}} + 1 \quad .
$$
This lead (in [GD]) to the {\it First Canonical Form}.
Another natural way is to use the {\it Fundamental Theorem of Arithmetic} and factor the integer into
prime powers, and then either write each prime as a sum of 1's and keep factorizing the exponents,
or write a prime as $1+(p-1)$ and factorize $p-1$ and continue recursively.
This lead, in [GD],  to the {\it Second Canonical Form}.

Either way, the {\it bottom line} is an {\it expression} that only uses $1$'s,  plus (``$+$''), times (``$*$''),
and exponentiation (``$\land$''). 

We will make the convention that $1$ can never be an argument of either multiplication or
exponentiation, or else there would be infinitely many ways of representing even $1$.

Given a positive integer $n$, how can we express it as a {\it formula}  {\bf only} using
the operations  $\{ +, *, \land \}$ and the integer $1$? 
[where we consider our operations as {\it binary}, i.e. {\it fan-in} $2$]

Of course there is only one way to express $1$, namely, $1$. There is also only one way to express $2$:
$$
2=1+1 \quad.
$$
[Strictly speaking we should write $2=(1)+(1)$, but we will abuse notation and abbreviate $(1)$ to $1$].

There are exactly two ways to express $3$
$$
3=(1+1)+1 \quad , \quad 3= 1+ (1+1) \quad .
$$
So far we only used addition. There are five ways to express $4$ only using addition:
$$
1+ ((1+1)+1) \quad , \quad 1+ (1+ (1+1)) \quad , \quad  (1+1)+(1+1) \quad , \quad ((1+1)+1)+1 \quad , \quad  (1+ (1+1) ) + 1 \quad .
$$

[In general there are $C_n=(2n)!/(n!(n+1)!)$ ways of expressing $n$ {\it only} using addition].

If you are also allowing multiplication, then we have, in addition (no pun intended)
$$
4=(1+1)*(1+1) \quad ,
$$
and if you are also allowing exponentiation, we have
$$
4=(1+1)\land (1+1) \quad .
$$

The above are examples of {\it formulas} whose inputs are always 1's.
The easiest way to define a formula is via  `grammars'.
If we only use addition, the additive formulas are given by the grammar
$$
F=1 \quad OR \quad (F)+(F)  \quad ,
$$
while the formulas that allow both addition and multiplication are defined by
$$
F=1 \quad OR \quad (F)+(F) \quad OR \quad (F)*(F) \quad ,
$$
and if you also allow exponentiation, then the grammar is
$$
F=1 \quad OR \quad (F)+(F) \quad OR \quad (F)*(F)  \quad OR \quad (F) \land (F) \quad .
$$

The above format is {\it infix}. As is well known (especially to  users of HP calculators)
one can get rid of parentheses, using {\it postfix} (alias Reverse Polish) notation.
The translation from infix to postfix is easy

$$
1 \rightarrow 1 \quad , \quad
a+b \rightarrow ab+ \quad , \quad
a*b \rightarrow ab* \quad , \quad
a \land b \rightarrow ab\land \quad .
$$
Of course these transformation rules are to be applied recursively.
For example, the expression
$$
(1+(1+1))\land((1+1)+1) \quad
$$
(representing twenty-seven), is written in postfix notation as
$$
111++11+1+\land \quad .
$$

We have already mentioned that the number of expressions of $n$ that only use addition, let's call it $C_a(n)$, 
is the famous Catalan sequence $(2n)!/(n!(n+1)!$ (why?). Let $C_{am}(n)$ be the number of such
expressions that use both addition and multiplication, and $C_{ame}(n)$ the number of expressions
that use the full arsenal of addition, multiplication, and {\bf exponentiation}.

In this short article (accompanied by a very long Maple package, and even longer sample output files) we
will answer the following questions.

$\bullet$ How to compute the sequences $C_{am}(n)$ and $C_{ame}(n)$ for as many $n$ as possible ?(it is unlikely that
there are  closed-form formulas).

$\bullet$ What is the asymptotics of $C_{am}(n)$ and $C_{ame}(n)$ as $n \rightarrow \infty$ ?

$\bullet$ How to draw {\it uniformly at random}, such an expression ?

$\bullet$ How to find the shortest possible expression for a given integer $n$. Of course, if you only use addition
all $C_a(n)$ expressions have the same length $2n-1$, but of course if one allows multiplication
one can get much shorter expressions, and if one also allows exponentiation, then one
can get yet shorter ones. [The length of such a minimal expression may be called the {\it computational complexity}
of the integer (w.r.t. the computational models discussed here)]

{\bf Enumeration}

{\bf Only using addition} 

Let $C_a(n)$ be the number of expressions for the positive integer $n$ {\it only} using
addition. Such an expression may be written as $n=k+(n-k)$ for some $1 \leq k < n$, and the number of these
is $C_a(k)C_a(n-k)$, so we have the non-linear recurrence
$$
C_a(n)= \sum_{k=1}^{n-1} C_a(k)C_a(n-k) \quad , \quad C_a(1)=1,
$$
whose solution is famously $(2n)!/(n!(n+1)!)$, the ubiquitous Catalan sequence [S] {\tt http://oeis.org/A000108}.

{\bf Using addition and multiplication}

Let $C_{am}(n)$ be the number of formula-trees with the leaves all $1$'s that represent the integer $n$, and
$C_{am}^{a}(n)$ be the number of those whose root is $+$, and $C_{am}^{m}(n)$ be the number of those whose root is $*$.
Then we have, of course
$$
C_{am}(n)=C_{am}^{a}(n)+C_{am}^{m}(n) \quad ,
$$
and the non-linear recurrences
$$
C_{am}^{a}(n) =\sum_{i=1}^{n-1} C_{am}(i)C_{am}(n-i) \quad ,
$$
$$
C_{am}^{m}(n) =\sum_{i>1, n/i \,\,\, integer}^{\lfloor n/2 \rfloor} C_{am}(i)C_{am}(n/i) \quad .
$$

[See procedures {\tt Cam(n)} and {\tt CamSeq(N)} in {\tt ArithFormulas}].

Using the {\it Wilf-methodology} [W][NW] we can use the {\it remembered} values of $C_{am}(n)$ and $C_{am}^{a}(n)$, $C_{am}^{m}(n)$
to generate {\it uniformly at random} such an expression. First use a {\it loaded} coin with probabilities
$C_{am}^{a}(n)/C_{am}(n)$, $C_{am}^{m}(n)/C_{am}(n)$ to decide whether the root-operation is ``plus'' or ``times'', and
in the former case use an $n-1$-faced loaded die whose faces are labeled $1, ..., n-1$, and the
probability of lending on $i$ is $C_{am}(i)C_{am}(n-i)/C_{am}^{a}(n)$, and continue recursively for 
$i,n-i$ assuming that it landed on $i$. Similarly if the loaded coin decided that the root-operation is ``$*$'' then,
create a loaded die whose faces are labeled by the non-trivial divisors of $n$, and the probability of
lending on face $i$ is $C_{am}(i)C_{am}(n/i)/C_{am}^{m}(n)$ and continue recursively.

[See procedures {\tt RaFamT(n)} and {\tt RaFamP(n)}in {\tt ArithFormulas}].

{\bf Using addition, multiplication and exponentiation}

Let $C_{ame}(n)$ be the number of formula-trees, whose internal nodes are in $\{+,*,\land \}$ and
whose leaves are all $1$'s, that represent the integer $n$, and
$C_{ame}^{a}(n)$ be the number of those whose root is $+$, $C_{ame}^{m}(n)$ be the number of those whose root is $*$,
$C_{ame}^{e}(n)$ be the number of those whose root is $\land$.

Then we have, of course
$$
C_{ame}(n)=C_{ame}^{a}(n)+C_{ame}^{m}(n)+C_{ame}^{e}(n) \quad ,
$$
and the non-linear recurrences
$$
C_{ame}^{a}(n) =\sum_{i=1}^{n-1} C_{ame}(i)C_{ame}(n-i) \quad ,
$$
$$
C_{ame}^{m}(n) =\sum_{i>1, n/i \,\,\, integer}^{\lfloor n/2 \rfloor} C_{ame}(i)C_{ame}(n/i) \quad .
$$
$$
C_{ame}^{e}(n) =\sum_{i^j=n, j>1} C_{ame}(i)C_{ame}(j) \quad .
$$

[See procedures {\tt Came(n)} and {\tt CameSeq(N)} in {\tt ArithFormulas}].

Using the {\it Wilf-methodology} [W][NW] we can use the {\it remembered} values of $C_{ame}(n)$ and 
$C_{ame}^{a}(n), C_{ame}^{m}(n), C_{ame}^{e}(n) $ to generate {\it uniformly at random} such an expression,
in an analogous way to the addition-multiplication trees above.

[See procedures {\tt RaFameT(n)} and {\tt RaFameP(n)}in {\tt ArithFormulas}].

{\bf Finding the Shortest Formula}

Using {\it Dynamical programming} we can find the shortest possible formula (measured in terms
of length in postfix notation), in either categories. We look at all the possible root operations and their subtrees
and pick the shortest possibility, using the previously obtained expressions for the children.

[See procedures {\tt ShortestTam(n)}, {\tt ShortestTame(n)} for the shortest formulas in infix (tree) notation and
procedures {\tt ShortestPam(n)}, {\tt ShortestPame(n)} for the shortest formulas in postfix (Reverse Polish) notation].

{\bf Asymptotics}

The well-known asymptotics for $C_a(n)=(2n)!/(n!(n+1)!)$ can be easily derived from Stirling's formula, yielding
$\frac{1}{\sqrt{\pi}} 4^n n^{-3/2}$. It is much harder to derive the asymptotics for
$C_{am}(n)$ and $C_{ame}(n)$ rigorously, but using procedure {\tt Zinn} of {\tt ArithFormulas}, we get
the following non-rigorous estimates

$$
C_{am}(n) \asymp c_1 n^{-3/2}(4.077...)^n \quad  ,
$$
$$
C_{ame}(n) \asymp c_2 n^{-3/2} (4.131...)^n \quad ,
$$
for some constants $c_1,c_2$.

{\bf The Book of Minimal Formulas}

To get the enumeration (up to $n=40$), and a list of optimal-length formulas for $n$ from $2$ to $8000$, generated by
procedures {\tt SeferAM(K1,K2)}  and {\tt SeferAME(K1,K2)} (with $K1=40,K2=8000$) for formulas using
only addition and multiplication and for formulas also using exponentiation, respectively, see the
two webbooks

{\tt http://www.math.rutgers.edu/\~{}zeilberg/tokhniot/oArithFormulas1} ,

{\tt http://www.math.rutgers.edu/\~{}zeilberg/tokhniot/oArithFormulas2} .

These minimal expressions are listed in postfix notation, ready to be entered into a Reverse Polish Calculator
(available on-line, e.g. {\tt http://www.alcula.com/calculators/rpn/}, viewed March 1, 2013).
They are given in the most memory-efficient way (using procedure {\tt MinMemory}) so as to minimize
the number of memory locations (stack-size) needed, i.e. realizing the {\it Strahler number}
(see {\tt Stra} in {\tt ArithFormulas}).

We also have analogous procedures for using {\it addition and exponentiation} (i.e. {\bf no} multiplication).
The output is presented in the following webbook

{\tt http://www.math.rutgers.edu/\~{}zeilberg/tokhniot/oArithFormulas3} .

\vfill\eject

{\bf Conclusion}

In addition to the great intrinsic interest of this project-what can be more natural or fundamental
than expressing integers?-it is also a {\it case study} in using {\bf Experimental Mathematics}
to enumerate, randomly generate, and optimally generate, combinatorial objects. We believe that
the same methodology could be applied to {\it Boolean formulas} and even {\it Boolean circuits},
that would shed yet another angle on the central problem of theoretical computer science, the notorious
P vs. NP problem. So far most of the work was done by humans, using pencil-and-paper. It is
about time that computers will put some effort towards settling the most central problem of their field,
or at the very least, give some {\it empirical} and {\it experimental} insight about it.

\bigskip

{\bf References}

[GD] Edinah K. Gnang and Patrick Devlin, { \it Some Integer Formula-Encodings and related algorithms},
to appear in Adv. Appl. Math. [available from {\tt www.arxiv.org}]

[I] Abraham Ibn Ezra, {\it ``Sefer HaEkhad''} [The ``Book of One''],  
available on-line from  \hfill\break
{\tt http://www.scribd.com/doc/17323522/Ibn-Ezra-Sefer-HaEchad} [viewed March 2, 2013]

[NW] Albert Nijenhuis and Herbert S. Wilf, {\it ``Combinatorial algorithms for computers and calculators''},
Academic Press, 2nd edition, 1978.

[S] Neil Sloane, {\it The On-Line Encyclopedia of Integer Sequences}, {\tt oeis.org}.

[W] Herbert S. Wilf,
{\it A unified setting for sequencing, ranking, and selection algorithms for combinatorial objects},
Advances in Mathematics {\bf 24} (1977), 281-291.

\bigskip

\hrule

\bigskip

Edinah K. Gnang, Computer Science Department, Rutgers University (New Brunswick), Piscataway, NJ 08854, USA.
{\tt gnang at cs dot rutgers dot edu}

Doron Zeilberger, Mathematics Department, Rutgers University (New Brunswick), Piscataway, NJ 08854, USA.
{\tt zeilberg at math dot rutgers dot edu}

{\bf March 4, 2013}

\end